\documentclass[3p,11pt]{elsarticle}
\usepackage{color}
\usepackage[colorlinks,linkcolor=blue,citecolor=blue]{hyperref}
\usepackage{amsmath,amsfonts,amssymb,amsthm}
\usepackage{mathrsfs} 
\usepackage{graphicx,epstopdf,float,subfigure,booktabs}
\usepackage[margin=10pt,font=small,labelfont=bf,labelsep=period]{caption}

\usepackage{multirow}
\usepackage{longtable}
\usepackage{mathtools}

\numberwithin{equation}{section} 

\biboptions{sort&compress}

\newtheoremstyle{mythm}{1.5ex plus 1ex minus .2ex}{1.5ex plus 1ex minus .2ex}
{\song}{\parindent}{\song\bfseries}{}{1em}{}
\theoremstyle{mythm}

\usepackage{algorithm}
\usepackage{listings,xcolor} 
\usepackage{textcomp}
\usepackage[framed,numbered,autolinebreaks,useliterate]{mcode} 

\DeclareCaptionFont{white}{\color{white}}
\DeclareCaptionFormat{listing}{\colorbox[cmyk]{0.43, 0.35, 0.35,0.01}{\parbox{\textwidth}{#1#2#3}}}
\newcommand{\bb}{\boldsymbol}

\begin{document}

\begin{frontmatter}

\title{Programming of linear virtual element methods in three dimensions}
\author{Yue Yu}
\ead{terenceyuyue@sjtu.edu.cn}
\address{School of Mathematical Sciences, Institute of Natural Sciences, MOE-LSC, Shanghai Jiao Tong University, Shanghai, 200240, P. R. China.}

\begin{abstract}
  We present a simple and efficient MATLAB implementation of the linear virtual element method for the three dimensional Poisson equation. The purpose of this software is primarily educational, to demonstrate how the key components of the method can be translated into code.
\end{abstract}


\end{frontmatter}


\section{Introduction}

The virtual element method (VEM), first introduced and analyzed in \cite{Beirao-Brezzi-Cangiani-2013}, is a generalization
of the standard finite element method on general polytopal meshes. In the past few years, people have witnessed
rapid progresses of virtual element methods (VEMs) for numerically solving various PDEs, see \cite{Ahmad-Alsaedi-Brezzi-2013,DeDios-Lipnikov-Manzini-2016,Zhao-Zhang-Chen-2018,Beirao-DaVeiga-Brezzi-2014} for
examples.
One can refer to \cite{Sutton-2017} for a transparent MATLAB implementation of the conforming
virtual element method for the Poisson equation in two dimensions.
The construction of the VEMs for three dimensional problems has been accomplished in many papers \cite{Gain-Talischi-Paulino-2014,Gain-Paulino-Leonardo-2015,Chi-Beirao-Paulino-2017,Beirao-Dassi-Russo-2017}. However, to the best of knowledge, no related implementation is publicly available in the literature.

In this paper, we are intended to present a simple and efficient MATLAB implementation of the linear virtual element method for the three dimensional Poisson equation on general polyhedral meshes. All important blocks are described in detail step by step. Although the current procedure is only for first-order virtual element spaces, the design ideas can be directly generalized to higher-order cases.

\section{Virtual element methods for the 3-D Poisson equation}

Let $\Omega \subset \mathbb{R}^3$ be a polyhedral domain and let $\Gamma$ denote a subset of its boundary consisting of some faces. We consider the following model problem
\begin{equation} \label{Poisson3D}
\begin{cases}
- \Delta u = r \quad & \mbox{in}~~\Omega, \\
u = g_D  \quad & \mbox{on} ~~\Gamma, \\
\partial_n u = g_N \quad & \mbox{on} ~~ \Gamma' = \partial \Omega \backslash \Gamma,
\end{cases}
\end{equation}
where $r\in L^2(\Omega)$ and $g_N\in L^2(\Gamma')$ are the applied load and Neumann boundary data, respectively, and $g_D\in H^{1/2}(\Gamma)$ is the Dirichlet boundary data function.

In what follows, we use $K$ to represent the generic polyhedral element with $f \subset \partial K$ being its generic face. The vertices of a face $f$ are in a counterclockwise order when viewed from the inside.
The virtual element method proposed in \cite{Beirao-Dassi-Russo-2017} for \eqref{Poisson3D} is to find $u_h \in V_\Gamma^k$ such that
\[a_h(u_h,v_h) = \ell_h(v_h), \quad v_h \in V_0^k,\]
where
\[a_h(u_h,v_h) = \sum\limits_{K\in \mathcal{T}_h} a_h^K(u_h,v_h), \quad
\ell_h(v_h) = \int_\Omega r_h v_h {\rm d}x + \int_{\Gamma'} g_h v_h {\rm d}s.\]
The local bilinear form is split into two parts:
\[a_h^K(v,w) = (\nabla v, \nabla w)_K + h_K S^K(v-\Pi_k^\nabla v, w - \Pi_k^\nabla w),\]
where $\Pi_k^\nabla : V^k(K) \to \mathbb{P}_k(K)$ is the standard elliptic projection, and $S^K$ is the stabilization term given as
\[S^K(v,w) = \sum\limits_{i=1}^{N^K} \chi_i(v) \chi_i(w),\]
where $\chi_i(v) = v(p_i)$ and $p_i$ is the $i$-th vertex of $K$ for $i=1,2,\cdots, N_K$. The local linear form of the right-hand side will be approximated as
\[\ell_h^K(v_h) = \int_K r \Pi_k^\nabla v_h {\rm d}x + \sum\limits_{f \subset \Gamma' \cap \partial K}\int_f g_N \Pi_{k,f}^\nabla v_h {\rm d}s,\]
where $\Pi_{k,f}^\nabla: V^k(f) \to \mathbb{P}_k(f)$ is the elliptic projector defined on the face $f$.

For the detailed introduction of the virtual element spaces, please refer to Section 2 in \cite{Beirao-Dassi-Russo-2017}. In this paper, we only consider the lowest order case $k=1$, but note that the hidden ideas can be directly generalized to higher order cases.

\section{Data structure and test script}

\begin{figure}[!htb]
  \centering
  \subfigure[Polyhedral mesh]{\includegraphics[height=6cm]{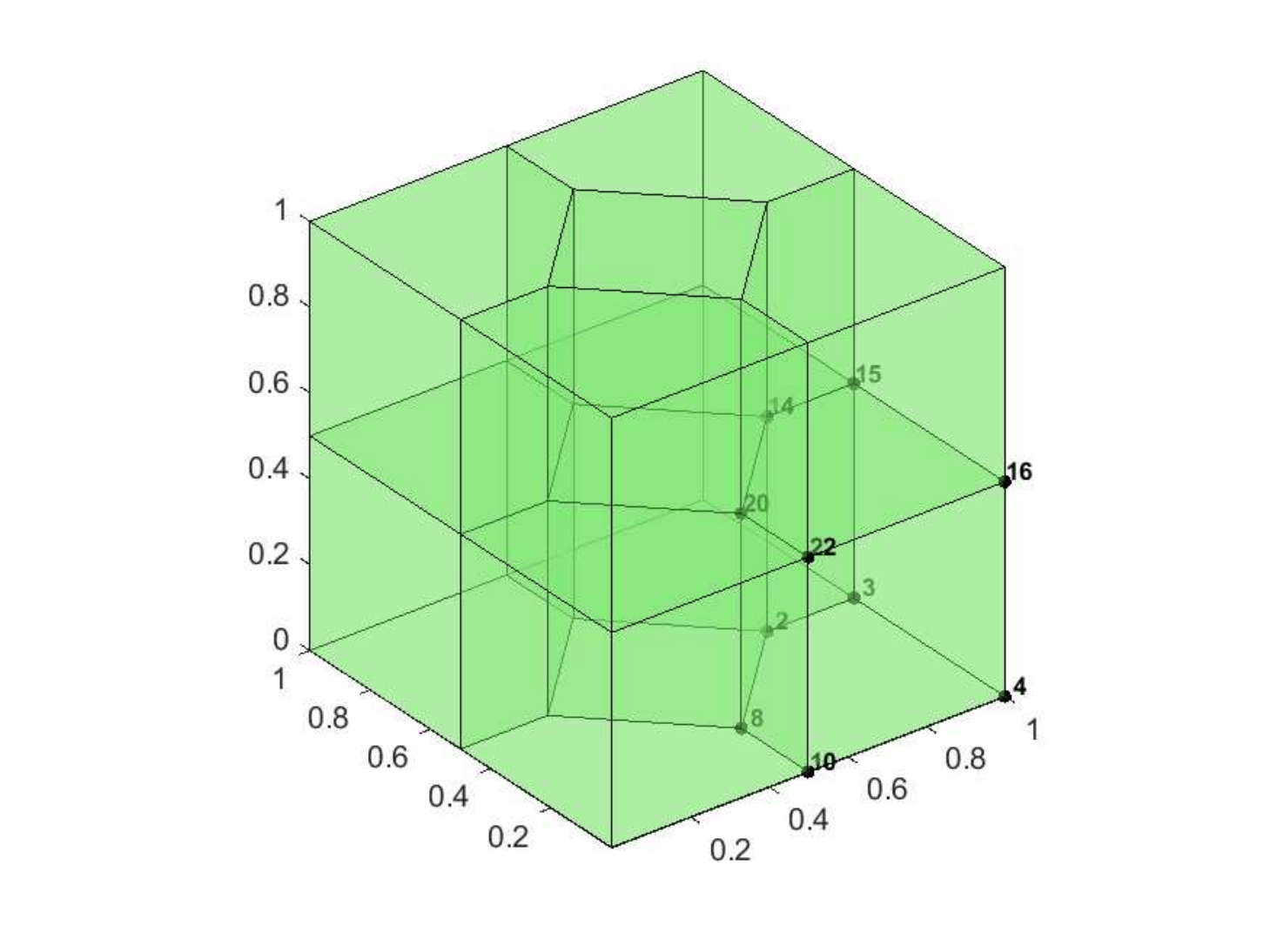}} \hspace{1cm}
  \subfigure[Representation of the first element]{\includegraphics[height=5cm]{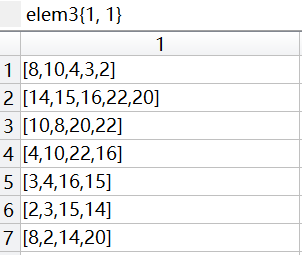}}\\
  \caption{Polyhedral mesh of a cube}\label{fig:mesh3ex1}
\end{figure}

We first discuss the data structure to represent polyhedral meshes. In the implementation, the mesh is represented by \mcode{node3} and \mcode{elem3}. The $\mcode{N} \times 3$ matrix \mcode{node3} stores the coordinates of all vertices in the mesh. \mcode{elem3} is a cell array with each entry storing the face connectivity, for example, the first entry \mcode{elemf = elem3\{1\}} for the mesh given in Fig.~\ref{fig:mesh3ex1}(a) is shown in Fig.~\ref{fig:mesh3ex1}(b), which is still represented by a cell array since the faces may have different numbers of vertices.

All faces including the repeated internal ones can be gathered in a cell array as

\vspace{0.5em}

\hspace{2cm} \mcode{allFace = vertcat(elem3\{:\});  \% cell }

\vspace{0.5em}

\noindent By padding the vacancies and using the \mcode{sort} and \mcode{unique} functions to rows, we obtain the face set \mcode{face}. The cell array \mcode{elem2face} then establishes the map of local index of faces in each polyhedron to its global index in face set \mcode{face}. The above structural information is summarized in the subroutine \mcode{auxstructure3.m}. The geometric quantities such as the diameter \mcode{diameter3}, the barycenter \mcode{centroid3} and the volume \mcode{volume} are computed by \mcode{auxgeometry3.m}. We remark that these two subroutines may be needed to add more information when dealing with higher order VEMs.

The test script is \mcode{main_PoissonVEM3.m} listed as follows. In the \mcode{for} loop, we first load the pre-defined mesh data, which immediately returns the matrix \mcode{node3} and the cell array \mcode{elem3} to the MATLAB workspace. Then we set up the Neumann boundary conditions to get the structural information of the boundary faces. The subroutine \mcode{PoissonVEM3.m} is the function file containing all source code to implement the 3-D VEM. When obtaining the numerical solutions, we calculate
the discrete $L^2$ errors and $H^1$ errors defined as
\[\mbox{ErrL2} = \sum\limits_{K\in \mathcal{T}_h} \| u - \Pi_1^\nabla u_h \|_{0,K}, \quad
  \mbox{ErrH1} = \sum\limits_{K\in \mathcal{T}_h} | u - \Pi_1^\nabla u_h |_{1,K}\]
by using respectively the subroutines \mcode{getL2error3.m} and \mcode{getH1error3.m}. The procedure is completed by verifying the rate of convergence through \mcode{showrateErr.m}.
\vspace{-0.5cm}
\begin{lstlisting}
%% Parameters
maxIt = 3;
h = zeros(maxIt,1);  N = zeros(maxIt,1);
ErrL2 = zeros(maxIt,1);
ErrH1 = zeros(maxIt,1);

%% PDE data
pde = Poisson3data();
bdNeumann = 'x==0'; % string for Neumann

%% Virtual element method
%[node3,elem3Refine,elem3] = cubemesh([0 1 0 1 0 1]);
for k = 1:maxIt
    % load mesh
    load( ['mesh3data', num2str(k), '.mat'] );
    %[node3,elem3Refine,elem3] = uniformrefine3(node3,elem3Refine);
    % get boundary information
    bdStruct = setboundary3(node3,elem3,bdNeumann);
    % solve
    [uh,info] = PoissonVEM3(node3,elem3,pde,bdStruct);
    % record
    N(k) = length(uh);  h(k) = (1/size(elem3,1))^(1/3);
    % compute errors in discrete L2, H1 and energy norms
    kOrder = 1;
    ErrL2(k) = getL2error3(node3,elem3,uh,info,pde,kOrder);
    ErrH1(k) = getH1error3(node3,elem3,uh,info,pde,kOrder);
end

%% Plot convergence rates and display error table
figure,
showrateErr(h,ErrL2,ErrH1);

fprintf('\n');
disp('Table: Error')
colname = {'#Dof','h','||u-u_h||','|u-u_h|_1'};
disptable(colname,N,[],h,'%0.3e',ErrL2,'%0.5e',ErrH1,'%0.5e');
\end{lstlisting}

In the following sections, we shall go into the details of the implementation of the 3-D VEM in \mcode{PoissonVEM3.m}.

\section{Elliptic projection on polygonal faces}

Let $f $ be a face of $K$ or a polygon embedded in $\mathbb{R}^3$. In the VEM computing, we have to get all elliptic projections $\Pi_{1,f}^\nabla \phi_f^T$ ready in advance, where $\phi_f$ is the nodal basis of the enhanced virtual element space $V^1(f)$ (see Subsection 2.3 in \cite{Beirao-Dassi-Russo-2017}).  To this end, it may be necessary to establish local coordinates $(s,t)$ on the face $f$.

\begin{figure}[!htb]
  \centering
  \includegraphics[scale=0.5]{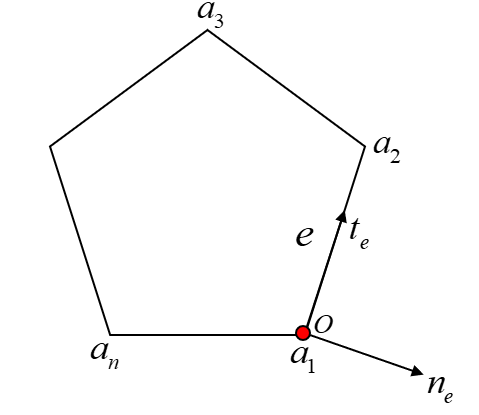}\\
  \caption{Local coordinate system of a face or polygon embedded in $\mathbb{R}^3$}\label{fig:localCord}
\end{figure}

As shown in Fig.~\ref{fig:localCord}, the boundary of the polygon is oriented in a counterclockwise order as $a_1, \cdots, a_n$. Let $e = a_1a_2$ be the first edge, and $n_e$ and $t_e$ be the normal vector and tangential vector, respectively. Then we can define a local coordinate system with $a_1$ being the original point by using these two vectors. Let $n_e = (n_1,n_2,n_3)$ and $t_e = (t_1,t_2, t_3)$. For any $a = (x,y,z) \in f$, its local coordinate $(s,t)$ is related by
\[\overrightarrow{Oa} = s\cdot n_e + t\cdot t_e, \quad \mbox{or} \quad
(x-x_1,y-y_1,z-z_1) = s \cdot (n_1,n_2,n_3) + t \cdot (t_1,t_2,t_3),\]
which gives
\[(s,t) = (x-x_1,y-y_1,z-z_1) \begin{bmatrix} n_1  &  n_2  & n_3 \\ t_1 & t_2  & t_2  \end{bmatrix}^{-1}, \]
with the inverse understood in the least squares sense.
When translating to the local coordinate system, we can compute all the matrices of elliptic projection in the same way for the Poisson equation in two-dimensional cases. For completeness, we briefly recall the implementation.
In what follows, we use the subscript ``$f$'' to indicate the locally defined symbols.

Let $\phi_f^T = [\phi_{f,1}, \cdots, \phi_{f,n}]$ be the basis functions of $V^1(f)$ and $m_f^T =[m_{f,1}, m_{f,2}, m_{f,3}]$ the scaled monomials on $f$ given as
\[m_{f,1} = 1, \quad m_{f,2} = \frac{s - s_f}{h_f}, \quad m_{f,3} = \frac{t - t_f}{h_f},\]
where $(s_f, t_f)$ and $h_f$ are the barycenter and the diameter of $f$, respectively. The vector form of the elliptic projector $\Pi_{1,f}^\nabla$ can be represented as
\begin{equation} \label{ellipticf}
\begin{cases}
(\nabla_f  m_f, \nabla_f \Pi_{1,f}^\nabla \phi_f^T)_f = (\nabla_f  m_f,  \nabla_f \phi_f^T )_f,  \\
P_0(\Pi_{1,f}^\nabla \phi_f^T) = P_0(\phi_f^T).
\end{cases}
\end{equation}
where
\[P_0(v) = \frac{1}{n}\sum\limits_{i=1}^n v (a_i).\]

Since $\mathbb{P}_1(f) \subset V^1(f)$, we can write
\[m_f^T = \phi_f^T \bb{D}_f, \qquad \bb{D}_f = (D_{i\alpha }), \quad D_{i\alpha } = \chi_{f,i} (m_{f,\alpha}),\]
where $\chi_{f,i}$ is the $i$-th d.o.f associated with $a_i$, and $\bb{D}_f$ is referred to as the transition matrix. We further introduce the following expansions
\[\Pi_{1,f}^\nabla{\phi_f^T} = \phi_f^T \bb{\Pi}_{1,f}^\nabla, \qquad
\Pi_{1,f}^\nabla{\phi_f^T} = m_f^T \bb{\Pi}_{1*,f}^\nabla.\]
One easily finds that
\[\bb{\Pi}_{1,f}^\nabla = \bb{D}_f\bb{\Pi}_{1*,f}^\nabla,\]
and \eqref{ellipticf} can be rewritten in matrix form as
\begin{align*}
\begin{cases}
\bb{G}_f \bb{\Pi}_{1*,f}^\nabla = \bb{B}_f,  \\
P_0(m_f^T) \bb{\Pi}_{1*,f}^\nabla = P_0(\phi_f^T)
\end{cases}, \quad \mbox{or denoted by} \quad
\tilde{\bb{G}}_f \bb{\Pi}_{1*,f}^\nabla = \tilde{\bb{B}}_f,
\end{align*}
where
\[\bb{G}_f = (\nabla_f  m_f, \nabla_f m_f^T)_f, \quad
\bb{B}_f = (\nabla_f  m_f, \nabla_f \phi_f^T)_f.\]
Note that the following consistency relation holds
\[\bb{G}_f = \bb{B}_f \bb{D}_f, \quad \tilde{\bb{G}}_f = \tilde{\bb{B}}_f \bb{D}_f.\]

Let \mcode{face} be the face set with internal faces repeated once. Then using the local coordinates we are able to derive all elliptic projections $\Pi_{1,f}^\nabla \phi_f^T$ as in 2-D cases. It is not recommended to carry out the calculation element by element in view of the repeated cost for the internal faces.

The above discussion is summarized in a subroutine with input and output as

\vspace{0.5em}

\hspace{2cm} \mcode{Pifs = faceEllipticProjection(P)},

\vspace{0.5em}

\noindent where \mcode{P} is the coordinates of the face $f$ and \mcode{Pifs} is the matrix representation $\bb{\Pi}_{1*,f}^\nabla$ of $\Pi_{1,f}^\nabla \phi_f^T$ in the basis $m_f^T$. One can derive all matrices by looping over the face set \mcode{face}:
\vspace{-0.5cm}
\begin{lstlisting}
%% Derive elliptic projections of all faces
faceProj = cell(NF,1);
for s = 1:NF
    % Pifs
    faces = face{s};  P = node3(faces,:);
    Pifs = faceEllipticProjection(P);
    % sort the columns
    [~,idx] = sort(faces);
    faceProj{s} = Pifs(:,idx);
end
\end{lstlisting}
Note that in the last step we sort the columns of $\bb{\Pi}_{1*,f}^\nabla$ in ascending order according to the numbers of the vertices. In this way we can easily find the correct correspondence on each element (see Lines 34-38 in the code of the next section).

The face integral is then given by
\begin{equation}\label{integralProj}
\int_f \Pi_{1,f}^\nabla \phi_f^T \mathrm{d}\sigma =  \int_f m_f^T \mathrm{d}\sigma \bb{\Pi}_{1*,f}^\nabla
=  ( |f|, 0, 0)\bb{\Pi}_{1*,f}^\nabla,
\end{equation}
where $|f|$ is the area of $f$ and the definition of the barycenter is used.

\section{Elliptic projection on polyhedral elements}

The 3-D scaled monomials $m^T = [m_1,m_2,m_3,m_4]$ are
\[m_1 = 1, \quad  m_2 = \frac{x-x_K}{h_K} , \quad  m_3 = \frac{y-y_K}{h_K} , \quad  m_4 = \frac{z-z_K}{h_K} ,\]
where $(x_K,y_K,z_K)$ is the centroid of $K$ and $h_K$ is the diameter, and the geometric quantities are computed by the subroutine \mcode{auxgeometry3.m}.  Similar to the 2-D case, we have the symbols $\bb{D}, \bb{G}, \tilde{\bb{G}}, \bb{B}$ and $ \tilde{\bb{B}}$. For example, the transition matrix is given by
\[\bb{D} = (D_{i \alpha}), \quad  D_{i \alpha} = \chi_i( m_\alpha) = m_\alpha (p_i).\]
The most involved step is to compute the matrix
\begin{align*}
\bb{B}
& = \int_K {\nabla m}  \cdot \nabla {\phi^T}\mathrm{d}x =  - \int_K {\Delta m}  \cdot {\phi^T}\mathrm{d}x + \sum\limits_{f \subset \partial K} {\int_f {(\nabla m \cdot {{\boldsymbol n}_f}){\phi^T}} } \mathrm{d}\sigma \\
& = \sum\limits_{f \subset \partial K} (\nabla m \cdot \bb{n}_f) \phi^T \mathrm{d}\sigma,
\end{align*}
where $\phi^T = [\phi_1,\phi_2,\cdots, \phi_{N_K}]$ are the basis functions with $\phi_i$ associated with the vertex $p_i$ of $K$.
According to the definition of $V^1(f)$, one has
\begin{align*}
\int_f (\nabla m \cdot \bb{n}_f) \phi^T \mathrm{d}\sigma
& = (\nabla m \cdot \bb{n}_f) \int_f \phi^T \mathrm{d}\sigma
  = (\nabla m \cdot \bb{n}_f) \int_f \Pi_{1,f}^\nabla \phi^T \mathrm{d}\sigma,
\end{align*}
and the last term is available from \eqref{integralProj}. Obviously, for the vertex $p_i$ away from the face $f$ there holds $\Pi_{1,f}^\nabla \phi_i = 0$. In the following code, \mcode{indexEdge} gives the row index in the face set \mcode{face} for each face of \mcode{elemf}, and \mcode{iel} is the index for looping over the elements.
\vspace{-0.5cm}
\begin{lstlisting}
    % ------- element information --------
    % faces
    elemf = elem3{iel};   indexFace = elem2face{iel};
    % global index of vertices
    index3 = unique(horzcat(elemf{:}));
    % local index of elemf
    Idx = min(index3):max(index3);
    Idx(index3) = 1:length(index3);
    elemfLocal = cellfun(@(face) Idx(face), elemf, 'UniformOutput', false);
    % centroid and diameter
    Nv = length(index3);  Ndof = Nv;
    V = node3(index3,:);
    xK = centroid3(iel,1); yK = centroid3(iel,2); zK = centroid3(iel,3);
    hK = diameter3(iel);
    x = V(:,1);  y = V(:,2);  z = V(:,3);

    % ------- scaled monomials ----------
    m1 = @(x,y,z) 1+0*x;
    m2 = @(x,y,z) (x-xK)/hK;
    m3 = @(x,y,z) (y-yK)/hK;
    m4 = @(x,y,z) (z-zK)/hK;
    m = @(x,y,z) [m1(x,y,z),m2(x,y,z),m3(x,y,z),m4(x,y,z)]; % m1,m2,m3,m4
    gradmMat = [0 0 0; 1/hK 0 0; 0 1/hK 0; 0 0 1/hK];

    % -------- transition matrix ----------
    D = m(x,y,z);

    % ----------- elliptic projection -------------
    B = zeros(4,Ndof);
    for s = 1:size(elemf,1)
        % --- information of current face
        % vertices of face
        faces = elemf{s};  P = node3(faces,:);
        % elliptic projection on the face
        idFace = indexFace(s);
        Pifs = faceProj{idFace}; % the order may be not correct
        [~,idx] = sort(faces);
        Pifs = Pifs(:,idx);
        % normal vector
        e1 = P(2,:)-P(1,:);  en = P(1,:)-P(end,:);
        nf = cross(e1,en); nf = nf./norm(nf);
        % area
        areaf = polyarea3(P);
        % --- integral of Pifs
        intFace = [areaf,0,0]*Pifs;  % local
        intProj = zeros(1,Ndof);  % global adjustment
        faceLocal = elemfLocal{s};
        intProj(faceLocal) = intFace;
        % add grad(m)*nf
        Bf = dot(gradmMat, repmat(nf,4,1),2)*intProj;
        B = B + Bf;
    end
    % constraint
    Bs = B;  Bs(1,:) = 1/Ndof;
    % consistency relation
    G = B*D;  Gs = Bs*D;
\end{lstlisting}

\section{Computation of the right hand side and assembly of the linear system}

The right-hand side is approximated as
\[F_K = \int_K f \Pi_1^\nabla \phi \mathrm{d} x = (\bb{\Pi}_{1*}^\nabla)^T \int_K f m \mathrm{d}x,\]
where $\Pi_1^\nabla$ is the elliptic projector on the element $K$ and $\bb{\Pi}_{1*}^\nabla$ is the matrix representation in the basis $m^T$. The integral $\int_K f m^T \mathrm{d}x$ can be approximated by
\[\int_K f m^T \mathrm{d}x = |K| f(\bb{x}_K) m(\bb{x}_K) = |K| f(\bb{x}_K) [1, 0, 0, 0]^T, \quad \bb{x}_K = (x_K, y_K, z_K). \]
One can also divide the element $K$ as a union of some tetrahedrons and compute the integral using the Gaussian rule. Please refer to the subroutine \mcode{integralPolyhedron.m} for illustration.

One easily finds that the stiffness matrix for the bilinear form is
\[\bb{A}_K = (\bb{\Pi}_{1*}^\nabla )^T \bb{G}\bb{\Pi}_{1*}^\nabla + h_K (\bb{I} - \bb{\Pi}_1^\nabla )^T(\bb{I} - \bb{\Pi}_1^\nabla ).\]
We compute the elliptic projections in the previous section and provide the assembly index element by element. Then the linear system can be assembled using the MATLAB function \mcode{sparse} as follows.
\vspace{-0.5cm}
\begin{lstlisting}
for iel = 1:NT

    ...

    % --------- local stiffness matrix ---------
    Pis = Gs\Bs;   Pi  = D*Pis;   I = eye(size(Pi));
    AK  = Pis'*G*Pis + hK*(I-Pi)'*(I-Pi);
    AK = reshape(AK,1,[]);  % straighten

    % --------- load vector -----------
    %fK = Pis'*[pde.f(centroid3(iel,:))*volume(iel);0;0;0];
    fun = @(x,y,z) repmat(pde.f([x,y,z]),1,4).*m(x,y,z);
    fK = integralPolyhedron(fun,3,node3,elemf);
    fK = Pis'*fK(:);

    % --------- assembly index for ellptic projection -----------
    indexDof = index3;
    ii(ia+1:ia+Ndof^2) = reshape(repmat(indexDof, Ndof, 1), [], 1);
    jj(ia+1:ia+Ndof^2) = repmat(indexDof(:), Ndof, 1);
    ss(ia+1:ia+Ndof^2) = AK(:);
    ia = ia + Ndof^2;

    % --------- assembly index for right hand side -----------
    elemb(ib+1:ib+Ndof) = indexDof(:);
    Fb(ib+1:ib+Ndof) = fK(:);
    ib = ib + Ndof;

    % --------- matrix for L2 and H1 error evaluation  ---------
    Ph{iel} = Pis;
    elem2dof{iel} = indexDof;
end
kk = sparse(ii,jj,ss,N,N);
ff = accumarray(elemb,Fb,[N 1]);
\end{lstlisting}

Note that we have stored the matrix representation $\bb{\Pi}_{1*}^\nabla$ and the assembly index \mcode{elem2dof} in the M-file so as to compute the discrete $L^2$ and $H^1$ errors.

\section{Treatment of the boundary conditions}

We first consider the Neumann boundary conditions. Let $f$ be a boundary face with $n$ vertices. The local load vector is
\[F_f = \int_f g_N \Pi_{1,f}^\nabla \phi_f \mathrm{d} \sigma
 = (\bb{\Pi}_{1*,f}^\nabla) \int_f g_N  m_f \mathrm{d} \sigma ,\]
where $g_N = \partial_{\boldsymbol n_f}u = \nabla u \cdot \bb{n}_f$. For simplicity, we provide the gradient $g_N = \nabla u$ in the PDE data instead and compute the true $g_N$ in the M-file. Note that the above integral can be transformed to a 2-D problem by using the local coordinate system as done in the following code, where \mcode{localPolygon3.m} realizes the transformation and returns some useful information, and \mcode{integralPolygon.m} calculates the integral on a 2-D polygon.
\vspace{-0.5cm}
\begin{lstlisting}
%% Assemble Neumann boundary conditions
bdFaceN = bdStruct.bdFaceN;  bdFaceIdxN = bdStruct.bdFaceIdxN;
if ~isempty(bdFaceN)
    Du = pde.Du;
    faceLen = cellfun('length',bdFaceN);
    nnz = sum(faceLen);
    elemb = zeros(nnz,1); FN = zeros(nnz,1);
    ib = 0;
    for s = 1:size(bdFaceN,1)
        % vertices of face
        faces = bdFaceN{s};   nv = length(faces);
        P = node3(faces,:);
        % elliptic projection on the face
        idFace = bdFaceIdxN(s);
        Pifs = faceProj{idFace}; % the order may be not correct
        [~,idx] = sort(faces);
        Pifs = Pifs(:,idx);
        % 3-D polygon -> 2-D polygon
        poly = localPolygon3(P);
        nodef = poly.nodef;  % local coordinates
        nf = poly.nf;   % outer normal vector of face
        centroidf = poly.centroidf;
        sc = centroidf(1); tc = centroidf(2);
        hf = poly.diameterf;
        Coord = poly.Coord; % (s,t) ---> (x,y,z)
        % g_N
        g_N = @(s,t) dot(Du(Coord(s,t)), nf);
        fun = @(s,t) g_N(s,t)*[1+0*s, (s-sc)/hf, (t-tc)/hf];
        Ff = integralPolygon(fun,3,nodef);
        Ff = Pifs'*Ff(:);
        % assembly index
        elemb(ib+1:ib+nv) = faces(:);
        FN(ib+1:ib+nv) = Ff(:);
        ib = ib + nv;
    end
    ff = ff + accumarray(elemb(:), FN(:),[N 1]);
end
\end{lstlisting}

The Dirichlet boundary conditions are easy to handle. The code is given as follows.
\vspace{-0.5cm}
\begin{lstlisting}
%% Apply Dirichlet boundary conditions
g_D = pde.g_D;  bdNodeIdxD = bdStruct.bdNodeIdxD;
isBdNode = false(N,1); isBdNode(bdNodeIdxD) = true;
bdDof = (isBdNode); freeDof = (~isBdNode);
nodeD = node3(bdDof,:);
u = zeros(N,1); uD = g_D(nodeD); u(bdDof) = uD(:);
ff = ff - kk*u;
\end{lstlisting}

In the above codes, \mcode{bdStruct} stores all necessary information of boundary faces. We finally derive the linear system
\mcode{kk*uh = ff}, where \mcode{kk} is the resulting coefficient matrix and \mcode{ff} is the right-hand side.
For small scale linear system, we directly solve it using the backslash command in MATLAB, while for large systems the algebraic multigrid method is used instead.
\vspace{-0.5cm}
\begin{lstlisting}
%% Set solver
solver = 'amg';
if N < 2e3, solver = 'direct'; end
% solve
switch solver
    case 'direct'
        u(freeDof) = kk(freeDof,freeDof)\ff(freeDof);
    case 'amg'
        option.solver = 'CG';
        u(freeDof) = amg(kk(freeDof,freeDof),ff(freeDof),option);
end
\end{lstlisting}
Here, the subroutine \mcode{amg.m} can be found in $i$FEM --- a MATLAB software package for the finite element methods \cite{iFEM}.

\section{Summary of the code}

\begin{figure}[!htb]
  \centering
  \subfigure[Polygonal mesh]{\includegraphics[scale=0.5]{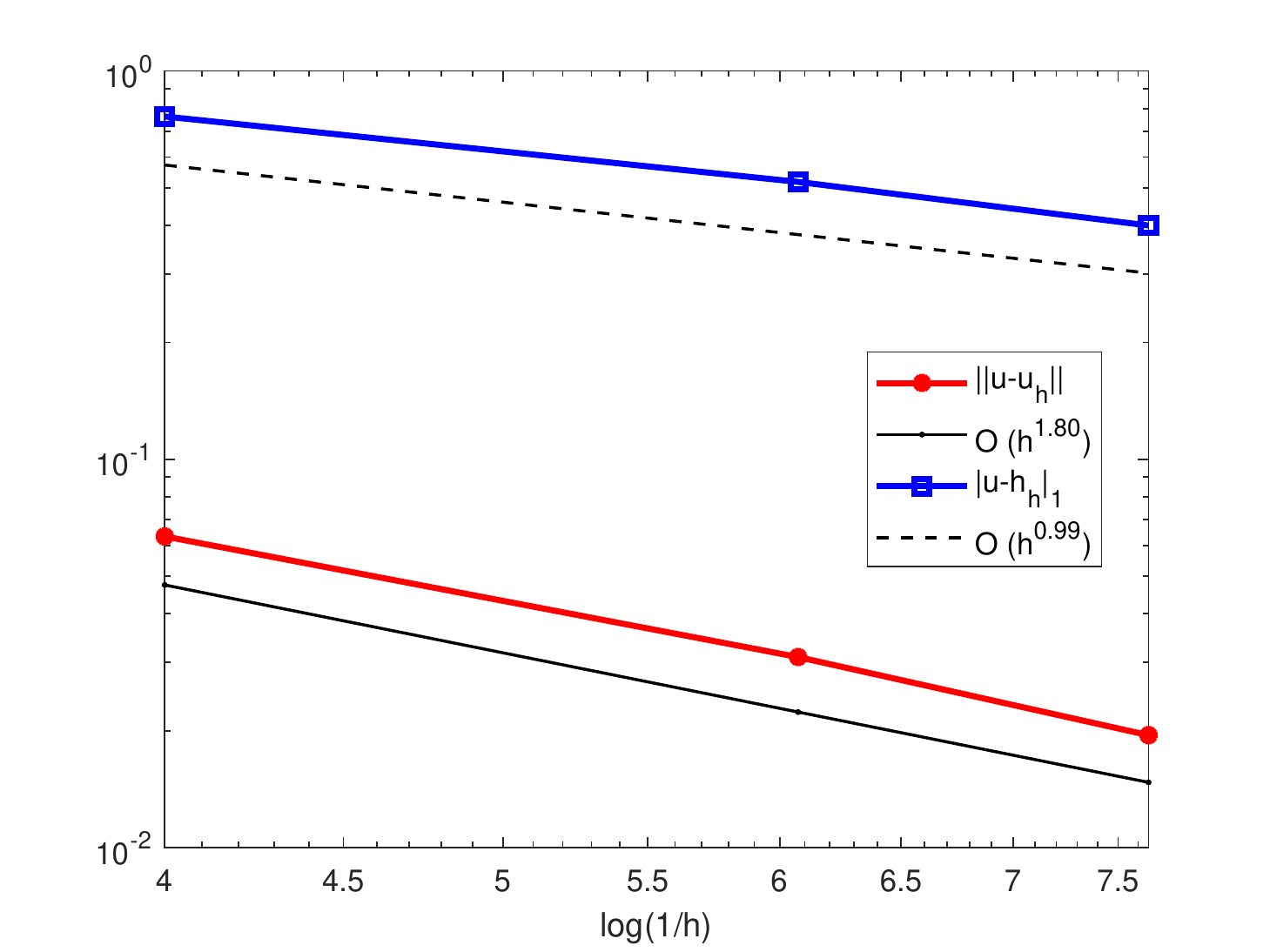}}
  \subfigure[Tetrahedral mesh]{\includegraphics[scale=0.5]{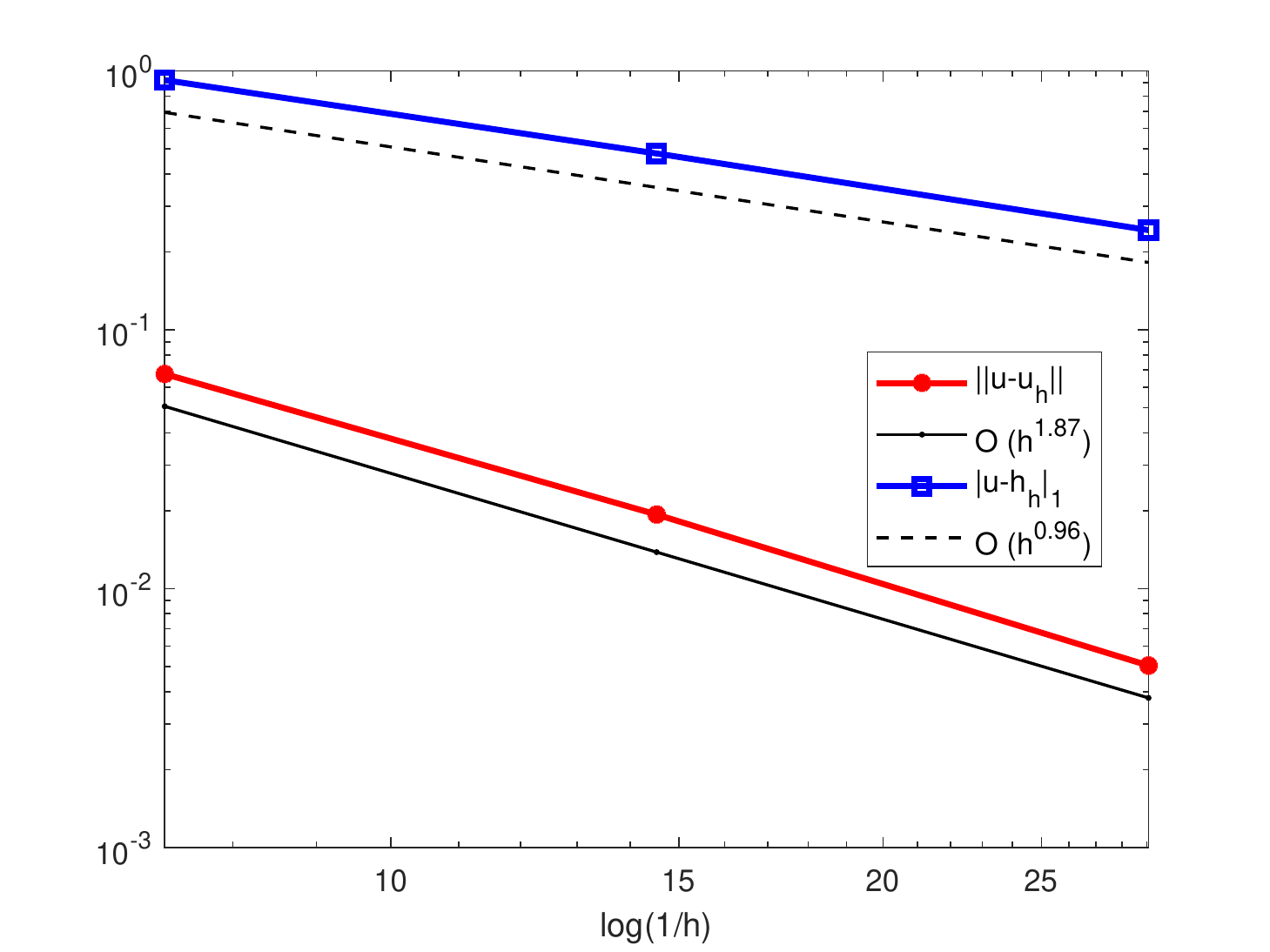}}\\
  \caption{Convergence rates in $L^2$ and $H^1$ norms for the 3-D Poisson equation}\label{fig:PoissonVEM3Rat}
\end{figure}

The complete M-file is \mcode{PoissonVEM3.m}. For the sake of the length, we only provide the main structure:
\vspace{-0.5cm}
\begin{lstlisting}
function [u,info] = PoissonVEM3(node3,elem3,pde,bdStruct)

%% Get auxiliary data

%% Derive elliptic projections of all faces

%% Compute and assemble the linear system

%% Assemble Neumann boundary conditions

%% Apply Dirichlet boundary conditions

%% Set solver

%% Store information for computing errors
\end{lstlisting}
To execute the M-file, one can refer to the example file \mcode{main_PoissonVEM3.m}.  We provide several simple meshes, such as \mcode{mesh3data1.mat} and \mcode{mesh3data2.mat}. All examples are implemented in MATLAB R2019b. The results for the test case 1 in \cite{Beirao-Dassi-Russo-2017} are displayed in Fig.~\ref{fig:PoissonVEM3Rat} and Tab.~\ref{tab:PoissonVEM3Err}, from which we observe the optimal rate of convergence both for the $H^1$ norm and $L^2$ norm.

\begin{table}[!htb]
  \centering
  \caption{Discrete $L^2$ and $H^1$ errors for polyhedral meshes}\label{tab:PoissonVEM3Err}
  \begin{tabular}{ccccccccccccccccc}
  \toprule[0.2mm]
   $N$            & $h$       &  ErrL2  & ErrH1 \\
  \midrule[0.3mm]
   170  & 2.500e-01  &   6.32804e-02   &  7.63490e-01 \\
   504  &  1.647e-01 &   3.09424e-02  &  5.18111e-01   \\
  1024  &  1.307e-01 &   1.94670e-02  &   4.00026e-01  \\
  \bottomrule[0.2mm]
\end{tabular}
\end{table}

Our code is available from GitHub (\url{https://github.com/Terenceyuyue/mVEM}) as part of the mVEM package which contains efficient and easy-following codes for various VEMs published in the literature. The test script {\color{gray} main\_PoissonVEM3.m} is stored in \mcode{vem3} folder for the current version.

\section*{Conflict of Interest}

The authors declare that they have no conflict of interest.











\end{document}